\documentclass[a4paper,english,12pt]{amsart}
\usepackage[all]{xy}
\usepackage{eucal}
\usepackage{mathrsfs}
\usepackage{enumerate}
\newtheorem{theorem}{ Theorem}[section]
\newtheorem{proposition}[theorem]{Proposition}
\newtheorem{lemma}[theorem]{ Lemma}
\newtheorem{corollary}[theorem]{ Corollary}
\newtheorem{definition}[theorem]{ Definition}
\theoremstyle{remark}

\newtheorem{example}[theorem]{Example}
\def \RM{\mathbb {R}}
\def \NM{\mathbb{N}}
\def \CM{\mathbb{C}}

\def \Bt {\mathcal{B}}

\def \Dt {\mathcal{D}}

\def \Ht {\mathcal{H}}

\def \abs{{\rm abs }\,}
\def \Id{{\rm Id }\,}
\def \d{\partial}

\def\dt{\delta} 
\def\a{\alpha}
\def\b{\beta}
\def\e{\varepsilon}

\def\G{\Gamma}   
\def \S{\Sigma}
\def \s{\sigma}

\def \to{\longrightarrow} 
\def \id{{\rm id\,}} 

\def \< {{\langle }}
\def \> {{\rangle }}

\newcommand{\Ot}{{{\mathcal O} }}

\begin{document}
\title [A generalisation of the Cauchy-Kovalevskaia theorem]{A generalisation of the\\ Cauchy-Kovalevskaia theorem} 
\author[Mauricio Garay]{Mauricio Garay}

\address{Mauricio Garay, \ Lyc\'ee Franco-Allemand\\
Rue Collin Mamet\\ 
78530 Buc.}
\email{garay91@gmail.com}
\date{Original version May 2012, (modified June 2015)}

\parskip=0.15cm


\begin{abstract}{ We prove that time evolution of a linear analytic initial value problem leads to sectorial holomorphic solutions in time.}
\end{abstract}
\maketitle
\section{Introduction} 
Among the class of systems of partial differential equations, evolutionary ones form a minority. They are nevertheless of considerable importance because they describe time evolution of physical data. This can already be seen for ordinary differential equations where, among differential equations, vector fields deserve a particular attention. The aim of this paper is to prove that evolutionary linear partial differential systems always admit sectorial analytic solutions, where the width of the sector depends on the regularity of the initial condition.

Before stating our main theorem, let us recall the results obtained by Kovalevska\"ia in her thesis~\cite{Kovalevskaia}~(see also \cite{Audin_Kovalevskaia} for historical aspects).
 
We consider the vector space $\CM^n$ with coordinates $z_1,\dots,z_n$ and let  $\Ot_n$ be the algebra of germs of holomorphic functions at  $0 \in \CM^n$ (series in $z_1,\dots,z_n$ which are analytic in some neighbourhood of the origin).
 We put
$$\d^I  :=\d_{z_1}^{i_1}\d_{z_2}^{i_2} \dots \d_{z_n}^{i_n} , \ I=(i_1,\dots,i_n) $$
and define the {\em order $\s(I)$ of the operator $\d^I$} as the sum of the coordinates of the vector $I \in \NM^n$:
$$\s(I):=i_1+\dots+i_n. $$
An {\em evolution equation} of order $s$ is a partial differential equation of the form
 $$\d_t x=F(x,\d_z^{I_1} x,\dots, \d_z^{I_k} x), \s(I_j) \leq s,\ x=(x_1,\dots,x_m),\ z=(z_1,\dots,z_n) $$
with some initial condition $x(t=0,\cdot)=x_0$, where $F$ is a polynomial mapping.

Kovalevska\"ia proved that the formal solution to such an initial value problem:
$$x(t,\cdot):=\sum_{k \geq 0} x_k t^k,\ x_k \in \Ot_n^m $$
exists and is unique. Then she proceeded to the analytic properties of time evolution. For $s=1$, she showed that the formal solution is holomorphic, in any sufficient small neighbourhood of the origin in $\CM^{n+1}$. A result now called the {\em Cauchy-Kovalevska\"{\i }a theorem}. For $s=2$, Kovalevska\"ia considered the particular case of the one dimensional {\em heat equation} and discovered that the formal solution might be  divergent.

To state {\em Kovalevska\"ia's heat equation theorem}, it is convenient to introduce the space $G^s_n$  of class $s$ Gevrey series in $n$ variables~\cite{Gevrey_CRAS,Gevrey_ANES}. These are formal power series: 
$$\sum_{I \in \NM^n} a_I z^I \in \CM[[z]] $$ 
such that
$$\sum_{I \in \NM^n} a_I \frac{z^I}{\left(\s(I)!\right)^{s-1}} $$
is convergent in a sufficiently small neighbourhood of the origin. For $s = 1$, Gevrey series are just analytic series, for $s<1$ these are entire functions and for $s>1$ these are divergent power series.
\begin{theorem}[\cite{Kovalevskaia}] The formal solution to the one dimensional heat equation
$$\d_t x=\d_{zz} x,\ x(t=0,-)=x_0,\ x_0 \in G_1^{s} $$
\begin{enumerate}[{\rm 1)}]
\item is a Gevrey class $2$ series;
\item has a unique holomorphic solution if $x_0  \in  G_1^{1/2}$;
\item is divergent if $x_0 \notin  G_1^{1/2}$.
\end{enumerate}
\end{theorem} 
Time evolution for the heat equation theorem became a classical subject and was treated in details in Hadamard's lectures on partial differential equations~\cite{Hadamard}.

In the eighties', Ouchi made an important step further, when he discovered that the divergent series associated to time evolution of a single linear partial differential equation are in fact asymptotic expansions of sectorial solutions~\cite{Ouchi}. We will extend the results of Kovalevska\"ia thesis and Ouchi's theorem to arbitrary systems of linear partial differential equations. 

\section{Statement of the theorem}
 According to the Kovalevska\"ia heat equation theorem, formal solutions will not be analytic in general. So we must look for {\em sectorial holomorphic solutions.} Let us clarify this notion.
 
  By {\em sector} of {\em width $\a \in [-\infty,+\infty]$}, we mean a subset of the form:
$$\S:=\{ z \in \CM^n :\ | z_i | \leq r ; | {\rm\ arg }\, z_i-\theta_i | \leq \frac{\a}{2} \} $$
for some {\em direction} $\theta=(\theta_1,\dots,\theta_n) \in (S^1)^n$ and some {\em radius} $r>0$. If the width is negative or infinite the condition on the angle is empty and the sector $\S$ is just a neighbourhood of the origin.

 We denote by $D_n(r) \subset \CM^n$ the polydisk:
$$D_n(r):=\{ z \in \CM^n: | z_i | \leq r,\ i=1,\dots,n \} $$
and by $D^*_n(r) \subset \CM$ the "pointed  polydisk":
$$D_n(r)^*:=\{ z \in (\CM^*)^n: | z_i | \leq r,\ i=1,\dots,n \}  $$

We denote by $G_n^\a(\theta)$ the algebra of functions which are holomorphic in an open subset containing some pointed polydisk of radius $r$ and which admit an asymptotic expansion of Gevrey class $\a$ inside any sector $\S$ of
width less than $\pi/(\a-1)$, direction $\theta$ and radius $r$.  

For instance the function $e^{-1/t}$ belongs to $G_1^2(0)$ which means that for any sector
$\S$ contained in the half plane 
$$\{ t \in \CM, {\rm Re}\, t >0 \} $$ the asympotic expansion of $e^{-1/t}$ is of Gevrey class $2$. This is indeed the case since it is equal to zero. Note that there is a map which associates to such a function its asymptotic expansion
$$G_n^\a(\theta) \to G_n^\a $$
and that it is not injective unless $\a \leq 1$.

Our theorem on time evolution requires some non-degeneracy condition.
We say that a linear partial differential operator of order $s$  is {\em non-degenerate}, if after a change of variables, it can be written in the form:
$$ \d_{z_n}^s -\sum_{\s(I) < s}A_I\d_z^I,\ A \in M(m,G_n^\a(\theta)). $$
Here $M(m,R)$ stands for $m \times m$ matrices with entries in the ring $R$ for some positive integer $m$.

When $\a=1$, for such a partial differential operator, the Cauchy-Kovalevska\"\i a theorem shows that there is a unique holomorphic solution with initial data
$$x(-,z_n=0)=x_0,\dots, \d_{z_n}^s x(-,z_n=0)=x_s. $$
For instance, an initial value problem of the form
$$\d_z^2 x=\d_t x,\ x(z=0,-)=x_0,\ \d_z x(z=0,-)=x_1$$ 
can be reduced to a system of order $1$~:
$$\left\{ \begin{matrix} \d_z x&=&y \\ 
\d_z y &=&\d_t x \end{matrix} \right.$$
with $x(t=0,-)=x_0,\ y(t=0,-)=x_1$. Therefore by Cauchy-Kovalevska\"ia theorem, it admits a unique holomorphic solution. This is  of minor interest for us, since we are interested in time evolution and not in spatial evolution. The main result of this paper is the:
\begin{theorem}
\label{T::main} Let $\sum_{\s(I) \leq s}A_I\d_z^I x$ be
a non-degenerate linear operator with  $A_I \in M(m,G_n^\a(\theta))$ for some direction $\theta$. The initial value problem
$$\d_t x=\sum_{\s(I) \leq s}A_I\d_z^I x,\ x_0:=x(t=0,-) \in \left( G_n^{\a}(\theta)\right)^m$$
admits solutions $x \in G_1^{\a s}(\theta') \hat  \otimes G_n^{\a}(\theta)$ for any $\theta' \in S^1$.
 \end{theorem}
 In the statement of the theorem we used a topological tensor product which can simply be understood as the space of functions in some pointed polydisk $D_{n+1}(r)^*$ which have asymptotic expansions of Gevrey class $\a s$ (resp. $\a$) in the $t$ variable (resp. $z$ variable) inside the sector $\S'$ (resp. $\S$). For more details on topological tensor products see~\cite{Grothendieck_PTT}.
\begin{example}
Consider the  Kovalevska\"ia example:
$$\d_t x=\d^2_z x,\ x(t=0,-)=\frac{1}{1-z}. $$
Here $\a=1,\ s=2$ so that $\S$ is of the form
$$\S=\S_1 \times D_2(r) $$
where $\S_1$ is a sector of width $\pi$, that is, a half plane. There exists a function
$$x(t,z):\{ (t,z) \in D_2(r) : {\rm Re}\, t>0\} \to \CM $$
which satisfies our initial value problem. We shall give examples of such functions in the next section. In real analysis, the solution of the heat equation in the circle can be solved by Fourier series and the flow is only defined for positive time. The situation is here completely different, since it admits solutions both for positive and negative time.
\end{example}

\begin{corollary}There is a unique holomorphic solution to an initial value problem of the form
$$\d_t x=\sum_{\s(I) \leq s}A_I\d_z^I x,\ x_0:=x(t=0,-) \in \left( G_n^{1/s}\right)^m,\ A_I \in M(m,G_n^{1/s}).$$
\end{corollary}
For the heat equation ($m=n=1$, $s=2$), we recover Kovalevska\"ia's theorem~(except for the statement on the divergence of the solution, which is discussed in the appendix).
The proof of Theorem~\ref{T::main} is based on a generalisation of Cauchy's {\em m\'ethode des majorantes} to general flows in infinite dimensional spaces.

\section{Formal evolution}
We start with the definition of formal evolution. In the linear case, this is quite obvious. Let $E$ be a vector space over a field $k$ and 
$$L:E \to E $$
a linear map. We denote by $E[[t]]:= k[[t]] \hat \otimes E$ the vector space of formal power series with coefficients in $E$:
$$\sum_{ n \geq 0} t^n \otimes a_n ,$$
that is the projective limit of the vector space $k[[t]]/(t^n) \otimes_k E  $.

The map $L$ induces a map $\id \otimes L$ on $E[[t]]$:
$$\left(\id \otimes L\right)(\sum_{ n \geq 0} t^n \otimes a_n)=\sum_{ n \geq 0} t^n \otimes L(a_n) .$$
We abuse notation and write $L$ for $\id \otimes L$. Similarly, we write $\d_t$ instead of $\d_t \otimes \Id$. We also identify $E$
with the subspace $1 \otimes E  \subset E[[t]]$.
If the field $k$ is of characteristic zero then the exponential is well-defined and $e^{t \otimes L}x_0$ is the unique solution to the initial value problem
$$\d_t x=L x,\ x(t=0,-)=x_0. $$
The operator $e^{tL}$ is called the {\em time evolution of the operator $L$.}

In order to extend this definition of evolution to non-linear operators, we first construct the Lie derivative. So let $E,F$ be locally convex vector spaces and let $U \subset E$ denote an open subset. 

 A mapping
$$f:E \supset U \to F, $$
is called {\em  G\^ateaux differentiable}   at a point $x \in U$, if for any $\xi \in E$, the following limits exists
$$D f(x)\xi:=\lim_{t \mapsto 0}\frac{f(x+t\xi)-f(x)}{t}. $$
Let
$$f:E \to F,\ X:E \to E$$ be  G\^ateaux differentiable mappings between locally convex spaces. The Lie derivative of $f$ along $X$ is defined by
$$ f \mapsto Df(x)X(x) .$$
We would like to find some vector space which is stable under the Lie derivative. In finite dimensional differential geometry, one may choose the space of $C^\infty$ function but these are difficult to handle in the infinite dimensional context, for general locally convex spaces. Therefore, we now assume that $k=\CM$ and consider the space of holomorphic maps from $E$ to $F$, denoted $\Ht(E,F)$. These are defined as follows:
 \begin{definition}
 A mapping $f:E \supset U \to F $  is called holomorphic if it satisfies the following two conditions:
\begin{enumerate}[{\rm i)}]
\item it is  continuous,
\item for any continuous linear mappings $j:\CM \to E$, $\pi:F \to \CM$ the map $\pi \circ f \circ j$ is holomorphic.
\end{enumerate}
\end{definition}
For instance, a linear mapping is holomorphic if and only if it is continuous.
Like in the finite dimensional case, holomorphic functions in infinitely many variables  are  convergent analytic power series (see \cite{Dineen} for more details).

 The elements of $\Ht(U,E)$ are called {\em holomorphic vector fields} in $U$.  
 By contracting the differential with a vector field $X \in \Ht(U,E)$, we define the {\em Lie derivative}
$$L_X:\Ht(U,E) \to \Ht(U,F),\ f \mapsto [x \mapsto Df(x)X(x)] $$
 for general locally convex spaces. As the Lie derivative is  linear map, we constructed in this way the derivation associated to a vector field
 in the infinite dimensional context.

The map 
$$L_X:\widetilde{E} \to  \widetilde{F},\ \widetilde{E}=\Ht(U,E),\  \widetilde{F}=\Ht(U,F)$$ 
being linear, we can consider the time evolution of any function. Now in the particular case $E=F$, we may consider the time evolution of the identity mapping. This defines in turn time evolution for general vector fields: 
\begin{definition} The formal flow of a holomorphic vector field $X \in \Ht(U,E)$ at $x_0$ is the evaluation  of the map $e^{tL_X}\Id$ at $x_0$, where $L_X$ is the Lie derivative along~$X$.
\end{definition}
Note that by construction the flow is a solution of the differential equation
$$\frac{dx}{dt}=(L_X \Id)(x)=X(x). $$
\begin{example}
Consider the {\em inviscid Burgers equation}:
$$\d_t x=x\d_z x,\ x(t=0,\cdot)=x_0 .$$
Denote by $\CM\{ z \}$ the vector space of convergent power series in one variable $z$, it has a natural topology~(see e.g. \cite{Grothendieck_EVT}). The vector field associated to our initial value problem is
$$\CM\{ z \} \to \CM\{ z \}, x \mapsto x\d_z x  .$$
The Lie derivative
$$L:\Ht(\CM\{ z \}) \to \Ht(\CM\{ z \}),\    [ x \mapsto f(x)] \mapsto [x \mapsto Df(x)x\d_z x]   $$
is linear and therefore admits a unique formal time evolution. Computation of time evolution up to order $2$ gives:
$$(e^{t \otimes L}\Id)x=x+t x\d_z x+\frac{t^2}{2}\left( 2x \left(\d_z x\right)^2+x^2\d_z^2 x \right) +o(t^2).$$
\end{example}
In simple words, the possibility to define differential calculus in the space of holomorphic functions $\Ht(\CM\{ z \})$ allows us to define formal flows like for the finite dimensional spaces. This explains the unicity of the formal solution to  initial value problems.
\section{Generalisation of the heat equation theorem}

We proceed to the Gevrey properties of formal solutions and first recall Borel resummation procedure~\cite{Borel_divergentes}~(see also~\cite{Balser_livre,Malgrange_resommabilite}).
Gevrey divergent series can be considered as asymptotic expansions of exponential integrals. For instance, the relation
$$\frac{1}{t}\int_0^{+\infty}e^{-\frac{\xi}{t}}\xi^kd\xi=k!t^{k} .$$
shows that for any polynomial, we get:
$$\int_0^{+\infty}e^{-\frac{\xi}{t}}\left( \sum_{k=0}^na_k\xi^k \right) d\xi=\sum_{k=0}^n k!a_kt^{k}. $$
More generally, if we consider a function
$$\xi \mapsto y(\xi)= \sum_{k\geq 0}a_k\xi^k $$
holomorphic in a neighbourhood of the real line which has at most exponential growth, then 
$$x(t)=\frac{1}{t}\int_0^{+\infty}e^{-\frac{\xi}{t}}y(\xi) d\xi$$
is a holomorphic function whose asymptotic expansion at the origin is:
$$\sum_{k \geq 0} k!a_kt^{k}. $$
The real line can of course be replaced by any closed curve which starts from $0$ and does not come back to it, for instance a real segment $[0,r]$. But along the real positive half-line, if the integrand is well-defined and decreases exponentially then the associated integral transformation maps the ring $\Dt_t$ of partial differential operator on $t$ to that on $\xi$ in the following way:
\begin{align*}\Dt_t &\to \Dt_\xi \\ \frac{1}{t} &\mapsto \d_{\xi} \\
\d_t &\mapsto \d_{\xi}  +\xi \d_\xi^2
\end{align*}
This procedure can be generalised to Gevrey classes $G^{1+s},\ s \geq  1$ using the exponential integral:
$$x(t):=\frac{1}{t^s}\int_0^{+\infty}e^{-\frac{\xi^s}{t^s}}y(\xi)d(\xi^s)$$

Let us now apply this resummation procedure to the heat equation
$$\d_t x=\d_{zz}x,$$
with Kovalevska\" \i a's initial value:
$$x_0:(\CM,0) \to (\CM,0),\ z \mapsto \frac{1}{1-z}.$$
The formal power series expansion of this Cauchy problem is of Gevrey class $2$
$$\hat x(t,z)=\frac{1}{1-z}\sum_{j,k}\frac{(2k)!}{k!}\frac{t^k}{(1-z)^{2k}}. $$
The substitution
$$k!t^{k}  \mapsto \xi^k $$
leads to the analytic series
 $$ y(\xi,z)=\frac{1}{1-z}\sum_{j,k}\frac{(2k)!}{(k!)^2}\frac{\xi^k}{(1-z)^{2k}}=\frac{1}{\sqrt{(1-z)^2-4\xi}}.$$
 The initial condition $x_0$ has a pole at $z=1$ and the integrand has itself a singularity at
$$\xi=\frac{(1-z)^2}{4}. $$ This is a particular case of the Lutz-Myiake-Sch\"afke theorem which states that the solution to the initial value problem "reproduces" the singularities of the initial data~\cite{Lutz}. Now let $\G$ be an arbitrary path starting at $0$, asymptotic to a line in the half plane
$$\{ \xi \in \CM: {\rm Re}\, \xi \geq 0 \}  $$
and which avoids the singularity $\xi=1/4$.
For $z$ sufficiently small, the power series $\hat x$ is the asymptotic expansion at $t=0$ of the function
 $$x_\G(z,t):=\frac{1}{t}\int_\G e^{-\frac{\xi}{t}}\frac{1}{\sqrt{(1-z)^2-4\xi}} d\xi.$$
 As explained above, the ring of partial differential operators $\Dt_{t,z}$ is mapped to $\Dt_{\xi,z}$. Thus our function
$$y(\xi,z)=\frac{1}{\sqrt{(1-z)^2-4\xi}} $$ 
is a solution of the partial differential equation
$$\left(\d_{\xi}  +\xi \d_\xi^2\right)y=\d^2_{z} y. $$
But 
$$\left(\d_t-\d^2_z \right)x_\G= \frac{1}{t}\int_\G\left(\d_{\xi}  +\xi \d_\xi^2-\d^2_{z}\right) y d\xi=0.$$
This means that  the functions $x_\G$ are solutions to our initial value problem, for any such choice of the path $\G$.  Note that the singularity
of the integrand implies that the solutions obtained by choosing different paths lead to different solutions and forces the divergence of the asymptotic series. This can already be observed for the case of the Euler equation:
$$t^2\frac{dx}{dt}+x=t,\ x(t)= \sum_{n \geq 0} (-1)^n n! t^{n+1}. $$
In this case the Borel transform is
$$\hat x(\xi)= \sum_{n \geq 0} (-1)^n \xi^{n+1}=\frac{1}{1+\xi} .$$
The singularity at $\xi=-1$ is responsible for the non unicity of the solution to our Cauchy problem~(see for instance \cite{Pham_livre_resurgence,Ecalle_fonctions}).

In any case, the first step to achieve Borel resummation of formal solutions is to ensure that they lie in some Gevrey class. This is established by the following 
\begin{theorem}
\label{T::Gevrey}  The formal solution to an initial value problem
$$\d_t x=\sum_{\s(I) \leq s}A_I\d_z^I x,\ x_0:=x(t=0,-) \in \left( G_n^{\a}\right)^m$$
of order $s$ is of Gevrey class $\a s$ in the time variable, that is, 
time evolution defines a map:
$$e^{tL}: \left( G_n^{\a}\right)^m  \to \left( G_n^{\a}\right)^m \hat \otimes G^{\a s}, x_0 \mapsto (e^{tL}\Id)x_0$$
where $L$ is the Lie derivative associated to the operator.
\end{theorem}
For a single partial differential equation ($m=1$), the theorem is again due to Ouchi. Using  techniques due to Boutet de Monvel and Kree, Yonemura simplified the proof~\cite{Boutet_Kree,Ouchi,Yonemura}.   Gevrey properties for some particular systems of partial differential equations (other than the heat equation) is proved in~\cite{Fernandez_Castro}. We postpone the proof of our theorem to the next section and first discuss the relation between formal solutions and asymptotic ones.

The  Borel transform
$$\Bt_s:G^s_1 \to G_1^1=\Ot_1,\ \sum_{k \geq 0} a_kt^{k}   \mapsto \sum_{k \geq 0} \frac{a_k}{(k!)^s}\xi^{k}$$
 associates a holomorphic function $y(\xi)$  to a Gevrey  series. In general this function is defined only in a small neighbourhood of the origin, therefore we cannot apply Borel resummation procedure to our formal solution as we did for the Kovalevska\"\i a example. So we replace, our integral formulas by
 $$x_r(t):=\frac{1}{t^s}\int_0^{r}e^{-\frac{\xi^s}{t^s}}y(\xi)d(\xi^s)$$
 where $r$ is now some finite number. The resulting function is still asymptotic to the formal solution $\hat x$. But under this new integral transformation, the differential operator $\d_t$ does not correspond to a differential operator in the $\xi$ variable anylonger. Therefore our function $x_r$ is, in general, not a solution to the initial value problem but rather a solution up to a flat function. For instance, for $r<1/4$, the asymptotic expansion of the function
 $$x_r(z,t):=\frac{1}{t}\int_0^r e^{-\frac{\xi}{t}}\frac{1}{\sqrt{(1-z)^2-4\xi}} d\xi$$
 is still a formal power series which satisfies the heat equation but the function itself does not.  In other words the function
 $$(\d_t-\d_{zz}) x_r $$
is flat at $t=0$ inside the half plane ${\rm Re}\, t>0$.
 
This might sound disappointing, and one might conclude that the above theorem gives simply no information about local analytic solutions, fortunately:
\begin{lemma} The existence of formal solutions implies the existence of asymptotic solutions, that is:\\
\begin{center}Theorem~\ref{T::Gevrey} $\implies$ Theorem~\ref{T::main}
\end{center}
 \end{lemma}
\begin{proof}
Let 
$$L:G_n(\theta)^m \to G_n(\theta)^m$$
be a linear partial differential operator of order $s \geq 2$. We consider an initial value problem
$$\d_t x=Lx,\ x(t=0,-)=x_0 \in G_n(\theta)^m. $$
Theorem~\ref{T::Gevrey} asserts that the formal solution is of Gevrey class $s$ in the time variable. It is therefore asymptotic expansion 
in any sector of width < $\pi/(s-1)$ of a holomorphic map $u$ defined in the pointed disk. The function $u$ is a solution to
our initial value problem up to a flat function. Substituting $x$ by $y+u $ in our system of partial differential equation, we get a new  system
$$ (*)\ \d_t y=L y+(Lu-\d_t u).$$
As the initial system is non-degenerate, in appropriate coordinates $L=\d_{z_n}^s-L'$ where $L'$ is of the form
$$L'= \sum_{\s(I)<s}A_I\d^I+\sum_{i=1}^{n-1}B_i \d_{z_i}^s.  $$
We re-write ($*$) as
 $$(**)\ \d_{z_n}^s y=\d_t y+L'y+g,\ g:=\d_t u-Lu.$$
 For any initial data:
 $$y_0=y(-,z_n=0),\ y_1=\d_{z_n}y(-,z_n=0),\dots, y_{s-1}=\d_{z_n}^{s-1}y(-,z_n=0),$$
 equation $(**)$ admits a unique holomorphic solution in the pointed disk. This can be seen, for instance, by applying the abstract Cauchy-Kovalevska\"\i a theorem to the Banach scale (see for instance~\cite{Baouendi,Nagumo,Nirenberg,Nishida,Ovsyannikov}):
 $$E_r:=C^0(R_r) \cap \Ot(R_r) $$
 with
 $$R_r:=\{ z \in \CM^n: r \leq z_i \leq 2r \} .$$

We  take the initial data $y_0=\dots=y_s=0$. The associated solution admits an asymptotic expansion equal zero when $t$ goes to zero inside the given sector. Indeed, by adding  $\d_{z_n}y ,\dots,\d_{z_n}^{s-1}y$ as new variables:
$$y_0=y,\ y_1=\d_{z_n}y,\dots,\ y_{s-1}=\d_{z_n}^{s-1}y. $$
 We reduce the system of partial differential equation to a first order system of the form
$$\d_{\tau} Y= TY+h ,\ Y=(y_0,\dots,y_{s-1}) $$
where $\tau=z_n$, $T$ is linear and $h$ is flat at $t=0$ in $\S$. The solutions 
are of the form
 $$Y=e^{\tau T}Y_0+\a,\ \a(\tau=0,-)=0$$
 where $\a$ satisfies the equation
 $$\d_{\tau} \a= T\a+h.$$
 As $h$ is flat, $\a=0$ is a formal solution to this equation, but the formal solution is unique thus $\a$ is flat at $t=0$ in $\S$. 
 As the initial condition is $Y_0=0$, the solution to our system of partial differential equations is constant $Y=\a$ and it is therefore flat. This shows that the holomorphic mapping 
$$x:=u+y$$ provides a solution to our original initial value problem in $G_n(\theta)^m$.  
\end{proof}
\section{Generalisation of Cauchy's majorant method}
 We consider vector fields in the  infinite dimensional Gevrey spaces $\left( G_n^s \right)^m$ and extend the classical Cauchy majorant method by comparing series in these different functional spaces. Via formal Borel transform these topological vector spaces are isomorphic to spaces of convergent power series. They are therefore endowed with a standard topology~(see~\cite{Grothendieck_EVT} for the definition of the topology).
 
Recall that a formal power series 
$$x:=\sum_{I \in \NM^n} a_I z^I \in \CM[[z]],\ z=(z_1,\dots,z_p)$$ is {\em majorated}
by another formal power series 
$$y:=\sum_{I \in \NM^n} b_Iz^I \in \RM_{+}[[z]]$$
if, for all $I \in \NM^n$, we have the estimates:
$$| a_I| \leq b_I.$$
In such cases, we use the notation
$$x \ll y.$$ 
In particular, $x \gg 0$ means that $x$ is a formal power series with real non-negative coefficients.
\begin{definition} Let $X,Y$ be two vector fields in $ \left(G_n^s \right)^m$.
A vector field $X$ in $  \left(G_n^s \right)^m $ majorates  another one $Y$ if  
$$x \gg y \implies X(x) \gg Y(y).$$
\end{definition}
In particular $X \gg 0$ means that:
$$x \gg 0 \implies X(x) \gg 0. $$
\begin{example}
Let $X$ be a vector field associated to a linear partial differential operator
$$X:(\Ot_n)^m \to (\Ot_n)^m,\ x \mapsto \sum_{\s(I) \leq s}A_I\d^I x. $$
Then $X\gg0$ provided that the entries of the matrices $A_I$ are analytic series with real positive coefficients, i.e., $A_I \gg 0$.
\end{example}
The following proposition is a direct consequence of the exponential formula for time evolution:
\begin{proposition}
\label{P::Taylor}
 Let $X,Y$ be two vector fields defined in an open subset of $ \left(G_n^s \right)^m $.
\begin{enumerate}[{\rm i)}]
\item If $X\ll Y$ then the flow of $X$ at $x_0 \gg 0$ is majorated by that of $Y$ at the same point,
\item If $X \gg 0$ and $y_0 \gg x_0$ then the flow of $X$ at $x_0$ is majorated by that of $X$ at $y_0$.
\end{enumerate}
\end{proposition}

We proceed to the proof of Theorem \ref{T::Gevrey} and start with a
\begin{proposition}\label{P::Gevrey}
 The following assertions are equivalent
\begin{enumerate}[{\rm 1)}]
\item the flow of any linear  initial value problem of order $s$ in $\left(G_n^\a\right)^m$ is of Gevrey class $\a s$ in the time variable ;
\item the flow of  any linear initial value problem of order $s$ in $\left(G_n^\a\right)^m$ is of Gevrey class $ \a s$ in the time variable ;
\item the flow of  any linear initial value problem  of order $s$ at $x_0=\sum_{n \geq 0}\left(n!\right)^{\a-1} z^n \in G_1^{\a}$ is of Gevrey class $s$ ;
\item  the flow of $x \mapsto x_0\d_z^sx  $ at the point $x_0=\sum_{n \geq 0}\left(n!\right)^{\a-1} z^n$ is of Gevrey class $\a s$.
\end{enumerate}
\end{proposition}
 \begin{proof}
Let us first make a remark. Consider a vector field defined by a linear differential operator:
$$X:\left( G_n^\a \right)^m \to \left( G_n^\a \right)^m ,\ x \mapsto \sum_{\s(I) \leq s} A_I\d_{z}^I x;\ I=(i_1,\dots,i_n)  $$
at a point $x_0$.
 
The map
$$ G_n^\a  \to G_n^\a,\ \a=\sum_I a_I z^I \mapsto \abs \a:=\sum_I |a_I| z^I$$
induces a map on matrices with coefficients in $G_n^\a$ that we denote in the same way.

Replace, in the initial value problem the $A_I$'s by $\abs A_I$'s and $x_0$ by $\abs x_0$.
By Proposition \ref{P::Taylor}, if the solution of this new initial value problem is of Gevrey class $k$
then $X,x_0$ has the same property. Therefore  it is sufficient to consider the case $X \gg 0$, $x_0 \in  \{x \gg 0 \}$.

\noindent $2) \implies 1)$.

Let us consider the linear mapping
$$\psi:\left( G_n^\a \right)^m \to G_n^\a,\ (x_1,\dots,x_m) \mapsto \sum_{k=1}^mx_k .$$
Write $A_I=(A_{I1},\dots,A_{Im}) \gg 0$ and put $f_I=\sum_k A_{Ik}$. For any $x\gg 0$, we have
$$\psi(\sum_{\s(I) \leq s} A_I\d_z^Ix)=\sum_{k=1}^m\sum_{\s(I) \leq s}A_{Ik}\d_z^Ix_k \ll \sum_{\s(I) \leq s} (\sum_{k=1}^mA_{Ik})\d_{z}^j(\sum_{k=1}^m x_k)=
\sum_{\s(I) \leq s} f_I \d_z^I \psi(x).$$
The exponential formula for time evolution implies that the image under $\psi$ of the flow of $X$  at $x_0$ is
majorated by the flow of $\sum_I f_I \d_z^I $ at $\psi(x_0)$.

\noindent $3) \implies 2)$.

Consider the open subset $U = \{ x \gg 0 \} \subset G_n^\a$.
The mapping
$$R:\CM \to \CM^n,\ z \mapsto (z,\dots,z) $$
induces a map 
$$R^*:G_n^\a[[t]] \supset U[[t]] \to G_1^\a[[t]],$$
and an element is of Gevrey class $ \b$ in the $t$ variable provided that it is the case of its image under $R^*$. 

The equalities
$$R^*\d_{z_i} z^k_j=k z^{k-1}\dt_{ij},\ \frac{d}{dz}R^* z_j^k=k z^{k-1} $$
give the estimate 
$$R^*\d_{z_i} \ll \frac{d}{dz}R^* .$$
Consider a vector field in $G_n^\a$ for the form: 
$$X:x \mapsto \sum_I f_I\d_{z}^Ix,\ f_I \gg 0$$
As  $R^*\d_{z_i} \ll \d_zR^* $, the flow of the vector field
$$G_1^\a \to G_1^\a,\ x \mapsto \sum_{\s(I) \leq s} R^*f_j \frac{d^{\s(I)} x}{d z^{\s(I)}}$$
at $R^*x_0, x_0 \gg 0$ majorates the image under $R^*$ of the flow of $X$.

Consider the Gevrey series
$$f(z):=\sum_{n \geq 0}\left(n!\right)^{\a-1} z^n \in G_1^{\a}.$$
and take
$$x_0=\sum_{n \geq 0} a_n z^n \in G_1^{\a}.$$
The series
$$ \sum_{n \geq 0} \frac{a_n}{(n!)^{\a-1}} z^n  $$
is analytic. Thus, by Hadamard's lemma, there exists $A,r>0$ such that
$$\frac{a_n}{(n!)^{\a-1}} \leq A r^n. $$ 
This means that the series $x_0$ is majorated by $Af(rz)$.  If $X \gg 0$, the formal flow passing through $x_0$ is majorated
by the formal flow passing through $Af(rz)$.  Up to multiplication of $x$ and $z$ by constants, we may assume that $A=r=1$.

\noindent $4) \implies 3)$.

Consider the flow of a vector field of the form:
$$X:x \mapsto \sum_{j=0}^s a_j\frac{d^jx}{dz^j}$$ at $f$ defined above.
 
As 
$$\frac{d^s}{dz^s}\gg \frac{d^j}{dz^j}$$
for any $j<s$, we get that the flow of $X$ at $x_0$
is majorated by that of the vector field
$$x \mapsto b(z)\frac{d^sx}{dz^s},\ b(z):=\sum_{j=0}^s a_j(z).$$
As before there exists constants $A,r>0$ such that 
$$ b(z) \ll Af(rz)  $$ 
and without loss of generality we may assume, as above,  that $A=r=1$. This concludes the proof of the proposition.
\end{proof}
To conclude the proof of Theorem \ref{T::Gevrey}, it remains to prove that the flow of the vector field
 $$X= f(z)\frac{d^s}{dz^s},\ f(z):=\sum_{n \geq 0}\left(n!\right)^{\a-1} z^n $$
 with initial condition
 $$x_0(z)=f(z)$$
is of Gevrey class $\a s$.

\begin{lemma} There exists a constant $C_\a>0$ such that for any real positive increasing sequence $(a_n)$ we have
$$fg \ll  C_\a  \sum_{n \geq 0} \left(n!\right)^{\a-1} a_n z^n $$
with
$$g:=\sum_{n \geq 0} \left(n!\right)^{\a-1}  a_n z^n .$$
\end{lemma}
\begin{proof}
Write
$$fg= \sum_{n \geq 0} c_n z_n$$
with
$$c_n:=\sum_{i+j=n} \left(i!\right)^{\a-1}\left(j!\right)^{\a-1}a_j \leq \left(\sum_{i+j=n} \left(i!\right)^{\a-1}\left(j!\right)^{\a-1}\right) a_n  .$$
One easily sees that
$$\sum_{i+j=n-1} i! j! \leq n! $$
therefore
$$\sum_{i+j=n} i! j! \leq 3 n! $$
For $\a\geq 1$, as $(a_n)$ is a positive increasing sequence, we get that:
$$c_n \leq  3^{\a-1}\left(n!\right)^{\a-1} a_n.$$
If $\a \leq 1$ we have
 $$\sum_{n \geq 0} \left(\sum_{i+j=n} \left(i!\right)^{\a-1}\left(j!\right)^{\a-1}\right) = (\sum_{i \geq 0} \left(i!\right)^{\a-1})^2<+\infty.$$
 Consequently the sequence
 $$ \left(\sum_{i+j=n} \left(i!\right)^{\a-1}\left(j!\right)^{\a-1}\right)_{n \in \NM} $$
 tends to zero at infinity. This implies that $c_n=o(a_n)$ and concludes the proof of the lemma.
\end{proof}
We have
$$\d_z^s f= \sum_{n \geq 0}\left(\frac{(n+s)!}{n!}\right)^\a\left(n!\right)^{\a-1} z^n .$$
The above lemma gives a constant $C_\a$ such that:
$$(f\d_z^s)^k f \ll  C_\a^k \sum_{n \geq 0}\left(\frac{(n+ks)!}{n!}\right)^\a\left(n!\right)^{\a-1} z^n.$$
Therefore
$$x(t)\ll  \sum_{n,k \geq 0}C_\a^k \left(\frac{(n+ks)!}{n!}\right)^\a\frac{\left(n!\right)^{\a-1}}{k!} z^nt^k.$$
Let us write $a_n \equiv b_n$ if the series $|a_n/b_n|$ is bounded by a geometric series. By Stirling's formula,  we have
$$\frac{(i+j)!}{i!j!} \equiv \frac{(i+j)^{i+j}}{i^i j^j}=e^{i\log (1+j/i)+j\log(1+i/j)} \leq e^{i+j} $$
thus
$$(i+j)! \equiv i!j! $$
We get that:
$$\left(\frac{(n+ks)!}{n!}\right)^\a\frac{\left(n!\right)^{\a-1}}{k!} \equiv \left(k!\right)^{\a s-1}\left(n!\right)^{\a-1} .$$
Thus the flow is of Gevrey class $\a s$ in the time variable.
 This concludes the proof of Theorem~\ref{T::Gevrey}. 
\appendix
\section{On the divergence of formal solutions}
 In~\cite{Lysik}, \L ysik proved a result similar to Kovalevska\"ia divergence result for the Korteweg--de Vries equation, namely that the solution
 to the initial value problem:
$$\d_t x=\d_{z}^{\, 3} x +x\d_z x,\ x(t=0,\cdot)= \frac{1}{1-z}  $$
is not holomorphic~(see also~\cite{Domrin}). More generally, one may wonder if our Gevrey estimate  for time evolution is optimal. This is indeed the case under very general asumptions:
\begin{theorem}
\label{T::divergence}
 Consider an evolutionary initial value problem of order~$s$
$$\d_t x=g(x,z)\d_{z_1}^sx+G(x,\d^{I_1} x,\d^{I_2} x,\dots,\d^{I_k} x ,\widehat{\d_{z_1}^sx}) ,\ \s( I_\a) \leq s, $$
with $x(t=0,\cdot)=x_0  $. Assume that $g,G \gg 0$ and $x_0 \gg 0$.
If the convergence radius of the formal Borel transform
$$\CM \to \CM^m,\ z_1 \mapsto B_{\a} x(z_1,0,\dots,0) $$
is finite then the formal solution to this initial value problem is not of Gevrey class $(\a s-\e)$, for any $\e>0$.
\end{theorem}  
\subsection*{Proof of Theorem \ref{T::divergence}}
The vector field associated to our initial value problem majorates
the vector field
$$X:x \mapsto g(x,z) \d_{z_1}^s $$
Moreover, the flow of  $X$ at  $x_0$ obviously majorates that of
$$Y:x \mapsto g(x_0,z)\d_z^j x$$
 at the same point. Finally, let 
$$\a z^I,\ \a \neq 0,\ I \in \NM^n$$ be a monomial appearing with a non-zero coefficient
in the Taylor expansion of $g$. We have
$$Y \gg \a z^I\d_{z_1}^j. $$
It remains to prove that the flow of
$$L=z^I  \d_{z_1}^s$$
at $x_0$ is not of Gevrey class $(\a s-\e)$, for any $\e>0$.

Given formal power series $f,g$, we write
$$f \succ g $$
if there are infinitely many coefficients of $f$ which are greater than that of $g$.
If $f \succ g$ and $g$ is not of Gevrey class $s$ then $f$ cannot be of Gevrey class $s$. 

Write 
$$x_0 = (x_{0,1},\dots, x_{0,m}).$$ 
Define
$$ f(z):=\sum_{n \geq 0}\left(n!\right)^{\a-1} z_1^n. $$
The assumption on $B_{\a s}$ implies that for at least one of the components of $x_0$, say
$x_{0,j}$,there exists $A,r>0$ such that : 
$$ x_{0,j} \succ  Af(rz_1) .$$
Up to a multiplication of $z_1$ and $x_0$ by constants, we may assume that $A=r=1$.

Now, the majorant
$$L^k \gg z^{kI}  \d_{z_1}^{ks},\ j \in \NM, $$

$$L^k f \gg z^{kI}  \sum_{n \geq 0}\left(\frac{(n+ks)!}{n!}\right)^\a\left(n!\right)^{\a-1} z_1^n \gg  z^{kI}  \sum_{n \geq 0}(k!)^{s\a}\left(n!\right)^{\a-1} z_1^n. $$
 Consequently
$$e^{tL_k}x_0 \succ \sum_{k \geq 0,n \geq 0} z^{kI} (k!)^{s\a-1}\left(n!\right)^{\a-1}  z_1^n t^k.$$
The right hand-side is not of Gevrey class $s\a-\e$ for any $\e>0$.
This proves the theorem.

\noindent {\em Acknowledgements.} I thank Boris Dubrovin for discussions from which this paper originated. Thanks also to Duco van Straten for encouragements and suggestions.
\bibliographystyle{amsplain}
\bibliography{master}

 \end{document}